\documentclass[11pt,letterpaper,reqno]{amsart}

\usepackage{amscd, amssymb, amsthm, amsmath}
\usepackage{hyperref} 
\usepackage{tikz}
\usepackage{ifthen}
\usetikzlibrary{positioning,calc}
\usepackage{mathrsfs}
\usepackage{ytableau}

\newtheorem{theorem}{Theorem}[section]

\newtheorem{lem}{Lemma}[section]

\newtheorem{exa}[lem]{Example}

\newtheorem{rem}{Remark}

\newtheorem*{sol}{Solution.}

\usepackage{mathtools}

\title[Revisiting the General Cubic: A Simplification of Cardano's Solution]{Revisiting the General Cubic: A Simplification of Cardano's Solution}

\author{Hsin-Chieh Liao}
\address{Hsin-Chieh Liao,~Department of Mathematics, University of Miami, Florida 33146, USA}
\email{h.liao@math.miami.edu}

\author{Mark Saul}
\address{Mark Saul,~Nine Nine Cultural and Educational Foundation, 9F, No. 300, Sec 3. Roosevelt Road, Zhongzheng Dist., Taipei City 10900, Taiwan}
\email{marksaul@earthlink.net}

\author{Peter J.-S Shiue}
\address{Peter J.-S Shiue,~Department of Mathematical Sciences,
University of Nevada, Las Vegas
Las Vegas, Nevada, 89154-4020, USA}
\email{shiue@unlv.nevada.edu}


\begin{document}




\maketitle

\setlength{\baselineskip}{15pt}

\setlength{\baselineskip}{15pt}

\section{Introduction}
Given a polynomial equation $x^3+ax^2+bx+c=0$  of degree 3 with real coefficients, we may translate the variable by replacing $x$ with $x-\frac{a}{3}$ to make the quadratic term vanish. We then obtain a simpler equation $x^3+px+q=0$ where $p=-\frac{a^2}{3}+b$ and $q=\frac{2a^3}{27}-\frac{ab}{3}+c$. Therefore, in order to solve a polynomial equation of degree 3, it is sufficient to solve for equations of the form $x^3+px+q=0$.\\

\noindent
In fact, the three roots of such a polynomial equations are
\[
	\sqrt[3]{A}+\sqrt[3]{B},~\omega\sqrt[3]{A}+\omega^2\sqrt[3]{B},~\omega^2\sqrt[3]{A}+\omega\sqrt[3]{B}
\]
where $\omega=\frac{-1+\sqrt{3}i}{2}$, $A=\frac{-q}{2}+\sqrt{\left(\frac{q}{2}\right)^2+\left(\frac{p}{3}\right)^3}$, $B=\frac{-q}{2}-\sqrt{\left(\frac{q}{2}\right)^2+\left(\frac{p}{3}\right)^3}$. This is the well-known formula found by the Italian mathematician G. Cardano(1501-76) \cite{HarrisStocker}. The formula appears in his book "Ars Magna". For the related historical details, we refer the readers to     
 \cite{DBurton}, \cite{Gilbert2015}.\\
 

Inspired by Sylvester's work \cite{JJSyv}, Chen in \cite{YCChen} proposed another way to solve the equation $x^3+px+q=0$. We set $p=-3rs$ and $q=rs(r+s)$.  Then the equation can be solved (if $r \neq s$) via the following new identity
\begin{equation}\label{chen}
	x^3-3rsx+rs(r+s)=\frac{s}{s-r}(x-r)^3+\frac{r}{r-s}(x-s)^3.
\end{equation}
Chen obtains the following results:
\begin{description}
\item[(i)] If $x^2-(r+s)x+rs=0$ has double roots $r=s$, the left hand side of (\ref{chen}) reduces to
\[
	x^3-3r^2x+2r^3=(x-r)^2(x+2r).
\]
\item[(ii)] If $r\neq s$ (where $s\neq 0$ since $p\neq 0$), then the three roots $x_1$, $x_2$, $x_3$ of $x^3+px+q=0$ are $x_i=\frac{r-su_i}{1-u_i}$ for $i=1,2,3$, where $u_1,u_2,u_3$ are the cube roots of $\frac{r}{s}$.
\end{description}

The goal of the paper is to continue Chen's work in \cite{YCChen} (see also Liao and Shiue\cite{HsinPeter}), studying in detail his method of solving the equation $x^3+px+q=0$.

\section{An alternative solution to the general cubic equation.}

\noindent
Although the case has been covered in \cite{YCChen}, we begin here, for completeness, with the case $r=s$, in which $x^3-3rsx+rs(r+s)=0$ can be factorized as follows:
\begin{align*}
	x^3-3r^2x+2r^3 &=(x^3-r^3)-3r^2(x-r)\\
				   &=(x-r)(x^2+rx-2r^2)\\
				   &=(x-r)^2(x+2r)
\end{align*}
Hence the roots are $r,r,-2r$, or equivalently $\sqrt{rs},\sqrt{rs},-2\sqrt{rs}$. Then we have the following theorem:

\begin{theorem}\label{r=s}
Let $x^3-3rsx+rs(r+s)=0$ be an equation with real coefficients. If $r=s$ then the three roots are
\[
	\sqrt{rs},\sqrt{rs},-2\sqrt{rs}.
\]
\end{theorem}

\begin{exa}
Solve $x^3-12x+16=0$.
\end{exa}
\begin{sol}

Because $rs=4$ and $rs(r+s)=16$, we have $r+s = 4$. Hence $r$ and $s$ are the two roots of $t^2-4t+4=0$. Solving this quadratic equation we obtain $r=s=2$. By theorem \ref{r=s} the three roots of $x^3-12x+16=0$ are $2,~2,~-4$.
\end{sol}







\begin{rem}
It is instructive to see how this relates to Cardano's formula.  The formula, as we have given it above, applies to the general cubic $x^3+px+q = 0$.  We work with the general form $x^3 -3rsx + rs(r+s)=0.$ So we have $p=-3rs, \ q = rs(r+s)$, or $rs = -\frac p 3, \ r+s= -\frac {3q} p$.  Hence $p$ and $q$ are roots of the quadratic $t^2 + \frac {3q}{4} t - \frac{p}{3} = 0$ . The discriminant of this quadratic equation is $\frac{9q^2}{p^2}+\frac{4p}{3}=\frac{1}{3p^2}(4p^3+27q^2)$. The term $4p^3+27q^2$ is exactly the classic discriminant of a cubic equation.
\end{rem}


\noindent
On the other hand, if $r\neq s$, from the identity (\ref{chen}) we get
\[
	\frac{s}{s-r}(x-r)^3+\frac{r}{r-s}(x-s)^3=x^3-3rsx+rs(r+s)=0,
\]
Then
\begin{equation}\label{Eq}
	s(x-r)^3=r(x-s)^3~~\Rightarrow~\left(\frac{x-r}{x-s}\right)^3=\frac{r}{s}.
\end{equation}
Since both $rs$ and $rs(r+s)$ are real, there are only two cases to be considered: either both $r$ and $s$ are real or $r$ and $s$ is a pair of complex conjugates.\\

\noindent
Suppose both $r$ and $s$ are real.  Then so is $\frac{r}{s}$. Hence we have $\frac{x-r}{x-s}=\sqrt[3]{\frac{r}{s}}, \sqrt[3]{\frac{r}{s}}\omega, \sqrt[3]{\frac{r}{s}}\omega^2$, where $\omega=e^{i\frac{2\pi}{3}}=\frac{-1+\sqrt{3}i}{2}$.\\
When $\frac{x-r}{x-s}=\sqrt[3]{\frac{r}{s}}$, then $\sqrt[3]{s}(x-r)=\sqrt[3]{r}(x-s)~\Rightarrow~ (\sqrt[3]{s}-\sqrt[3]{r})x=\sqrt[3]{s}r-\sqrt[3]{r}s$.  Therefore we get
\begin{align*}
	x=\frac{\sqrt[3]{s}r-\sqrt[3]{r}s}{\sqrt[3]{s}-\sqrt[3]{r}}=\frac{r^{\frac{1}{3}}s^{\frac{1}{3}}(r^{\frac{2}{3}}-s^{\frac{2}{3}})}{s^{\frac{1}{3}}-r^{\frac{1}{3}}}=-r^{\frac{1}{3}}s^{\frac{1}{3}}(r^{\frac{1}{3}}+s^{\frac{1}{3}}).
\end{align*}
Similarly, when $\frac{x-r}{x-s}=\sqrt[3]{\frac{r}{s}}\omega$, we have $\sqrt[3]{s}(x-r)=\sqrt[3]{r}\omega(x-s)~\Rightarrow~(\sqrt[3]{s}-\omega\sqrt[3]{r})x=\sqrt[3]{s}r-\omega\sqrt[3]{r}s$.  It follows that 
\begin{align*}
	x &=\frac{\sqrt[3]{s}r-\omega\sqrt[3]{r}s}{\sqrt[3]{s}-\omega\sqrt[3]{r}}=\frac{r^{\frac{1}{3}}s^{\frac{1}{3}}(r^{\frac{2}{3}}-\omega s^{\frac{2}{3}})}{s^{\frac{1}{3}}-\omega r^{\frac{1}{3}}}=
	\frac{r^{\frac{1}{3}}s^{\frac{1}{3}}(r^{\frac{2}{3}}-\omega s^{\frac{2}{3}})\omega^2}{\left(s^{\frac{1}{3}}-\omega r^{\frac{1}{3}}\right)\omega^2}\\
	&=-\frac{r^{\frac{1}{3}}s^{\frac{1}{3}}(\omega^2 r^{\frac{2}{3}}- s^{\frac{2}{3}})}{\omega^2\left(\omega r^{\frac{1}{3}}-s^{\frac{1}{3}}\right)}=-\omega r^{\frac{1}{3}}s^{\frac{1}{3}}(\omega r^{\frac{1}{3}}+s^{\frac{1}{3}})\\
	&=-r^{\frac{1}{3}}s^{\frac{1}{3}}(\omega^2 r^{\frac{1}{3}}+\omega s^{\frac{1}{3}}).
\end{align*}
When $\frac{x-r}{x-s}=\sqrt[3]{\frac{r}{s}}\omega^2$,  a similar argument yields $x=-r^{\frac{1}{3}}s^{\frac{1}{3}}(\omega r^{\frac{1}{3}}+\omega^2 s^{\frac{1}{3}})$. Summarizing the discussion above we obtain the following theorem:

\begin{theorem}\label{Case1}
Let $x^3-3rsx+rs(r+s)=0$ have real coefficients with $r\neq s$. If both $r$ and $s$ are real, then the three roots are
\[
	-r^{\frac{1}{3}}s^{\frac{1}{3}}(r^{\frac{1}{3}}+s^{\frac{1}{3}}),
	~-r^{\frac{1}{3}}s^{\frac{1}{3}}(\omega^2 r^{\frac{1}{3}}+\omega s^{\frac{1}{3}}),
	~-r^{\frac{1}{3}}s^{\frac{1}{3}}(\omega r^{\frac{1}{3}}+\omega^2 s^{\frac{1}{3}}).
\]
The equation has one real root and a pair of complex conjugate roots.
\end{theorem}

\begin{exa}
Find all solutions to $x^3-6x-9=0$.
\end{exa}
\begin{sol}
The numbers $r$ and $s$ are the roots of the quadratic equation $t^2+\frac{9}{2}t+2=0$.  Solving this  equation we get $r=-\frac{1}{2}$ and $s=-4$ which are two distinct real numbers.  Hence by theorem \ref{Case1} the three solutions to the original equations are
\[
	-2^{\frac{1}{3}}\left(\left(-\frac{1}{2}\right)^{\frac{1}{3}}+(-4)^{\frac{1}{3}}\right)=-((-1)^{\frac{1}{3}}+(-8)^{\frac{1}{3}})=3,
\]
\[
	-2^{\frac{1}{3}}\left(\left(\frac{-1-\sqrt{3}i}{2}\right)\left(-\frac{1}{2}\right)^{\frac{1}{3}}+\left(\frac{-1+\sqrt{3}i}{2}\right)(-4)^{\frac{1}{3}}\right)=-\left(\frac{1+\sqrt{3}i}{2}+\frac{2-2\sqrt{3}i}{2}\right)=\frac{-3+\sqrt{3}i}{2},
\]
\[
	-2^{\frac{1}{3}}\left(\left(\frac{-1+\sqrt{3}i}{2}\right)\left(-\frac{1}{2}\right)^{\frac{1}{3}}+\left(\frac{-1-\sqrt{3}i}{2}\right)(-4)^{\frac{1}{3}}\right)=-\left(\frac{1-\sqrt{3}i}{2}+\frac{2+2\sqrt{3}i}{2}\right)=\frac{-3-\sqrt{3}i}{2}.
\]
\end{sol}





\noindent
Next let us consider the second case, when $r$, $s$ is a pair of complex conjugates. We first establish some notation.
\ 
\ 

Any complex number $z=x+iy\neq 0$ can be expressed as $z=|z|e^{i\theta}=|z|(\cos\theta+i\sin\theta)$, where $\theta={\rm Arg}(z)$ is the argument of $z$ in the interval $(-\pi,\pi]$. We write $\overline{z}=x-iy$ for the complex conjugate of $z$ and ${\rm Re}(z)$ for the real part of $z$.
\ 
\ 

For any positive integer $m$, the complex number $z$ has $m$ $m$-th roots $\alpha,\alpha\zeta,\alpha\zeta^2,\ldots,\alpha\zeta^{m-1}$, where $\alpha=\sqrt[m]{|z|}e^{i\frac{\theta}{m}}$ and $\zeta=e^{i\frac{2\pi}{m}}$. In this paper we write $\alpha$ as $z^{\frac{1}{m}}$ . In particular, when $\theta=0$, $z^{\frac{1}{m}}=\sqrt[m]{z}\in\mathbb{R}$.\\

From (\ref{Eq}), we get $\frac{x-r}{x-s}=\left(\frac{r}{s}\right)^{\frac{1}{3}},~\left(\frac{r}{s}\right)^{\frac{1}{3}}\omega,~\left(\frac{r}{s}\right)^{\frac{1}{3}}\omega^2$. After a similar computation, we see that the solutions are given by the same expressions as in the former case:
\begin{equation}\label{Roots}
	-r^{\frac{1}{3}}s^{\frac{1}{3}}(r^{\frac{1}{3}}+s^{\frac{1}{3}}),
	~-r^{\frac{1}{3}}s^{\frac{1}{3}}(\omega r^{\frac{1}{3}}+\omega^2 s^{\frac{1}{3}}),
	~-r^{\frac{1}{3}}s^{\frac{1}{3}}(\omega^2 r^{\frac{1}{3}}+\omega s^{\frac{1}{3}}).
\end{equation}
However, when $r,s$ are complex conjugates, the expressions can be further simplified as follows: 

\ \ 
Let $r=|r|e^{i\theta}$ and $s=|s|e^{-i\theta}$ with $|r|=|s|$.  Then $r^{\frac{1}{3}}=|r|^{\frac{1}{3}}e^{i\frac{\theta}{3}}=(rs)^{\frac{1}{6}}e^{i\frac{\theta}{3}}$,  $s^{\frac{1}{3}}=(rs)^{\frac{1}{6}}e^{-i\frac{\theta}{3}}=\overline{r^{\frac{1}{3}}}$, and we can rewrite the solutions as
\[
	-r^{\frac{1}{3}}s^{\frac{1}{3}}(r^{\frac{1}{3}}+s^{\frac{1}{3}})=-(rs)^{\frac{1}{3}+\frac{1}{6}}\left(e^{i\frac{\theta}{3}}+e^{-i\frac{\theta}{3}}\right)=-2\sqrt{rs}\cos\frac{\theta}{3},
\]
\[
	-r^{\frac{1}{3}}s^{\frac{1}{3}}(\omega r^{\frac{1}{3}}+\omega^2 s^{\frac{1}{3}})
	=-(rs)^{\frac{1}{3}}(\omega r^{\frac{1}{3}}+\overline{\omega r^{\frac{1}{3}}})=-2(rs)^{\frac{1}{3}}{\rm Re}(\omega r^{\frac{1}{3}})=-2\sqrt{rs}\cos\left(\frac{\theta}{3}+\frac{2\pi}{3}\right),
\]
\[
	-r^{\frac{1}{3}}s^{\frac{1}{3}}(\omega^2 r^{\frac{1}{3}}+\omega s^{\frac{1}{3}})
	=-(rs)^{\frac{1}{3}}(\omega^2 r^{\frac{1}{3}}+\overline{\omega^2 r^{\frac{1}{3}}})=-2(rs)^{\frac{1}{3}}{\rm Re}(\omega^2 r^{\frac{1}{3}})=-2\sqrt{rs}\cos\left(\frac{\theta}{3}+\frac{4\pi}{3}\right).
\]

\begin{theorem}\label{Case2}
Let $x^3-3rsx+rs(r+s)=0$ have real coefficients with $r\neq s$. If $r$ and $s$ are a pair of complex conjugates and we let $\theta={\rm Arg}(r)$,  then the three roots of the equation are
\[
	-2\sqrt{rs}\cos\frac{\theta}{3},~-2\sqrt{rs}\cos\left(\frac{\theta}{3}+\frac{2\pi}{3}\right),~-2\sqrt{rs}\cos\left(\frac{\theta}{3}+\frac{4\pi}{3}\right).
\]
The equation has three real roots.
\end{theorem}

\begin{rem}\label{root}
Note that from the result of theorem \ref{Case2}  we can derive the following trigonometric identities by Vieta's formula,
\[
	\cos\frac{\theta}{3}+\cos\left(\frac{\theta}{3}+\frac{2\pi}{3}\right)+\cos\left(\frac{\theta}{3}+\frac{4\pi}{3}\right)=0
\]
\[
	\cos\frac{\theta}{3}\cos\left(\frac{\theta}{3}+\frac{2\pi}{3}\right)+\cos\frac{\theta}{3}\cos\left(\frac{\theta}{3}+\frac{4\pi}{3}\right)+\cos\left(\frac{\theta}{3}+\frac{2\pi}{3}\right)\cos\left(\frac{\theta}{3}+\frac{4\pi}{3}\right)=-\frac{3}{4}
\]
\[
	\cos\frac{\theta}{3}\cos\left(\frac{\theta}{3}+\frac{2\pi}{3}\right)\cos\left(\frac{\theta}{3}+\frac{4\pi}{3}\right)=-\frac{\cos\theta}{4}.
\]
The last identity follows from the fact that the product of the three roots is equal to $-rs(r+s)$ and $r=|r|e^{i\theta}, s=|r|e^{-i\theta}$.
\end{rem}

\begin{rem}
In remark \ref{root},  we have $\cos\frac{\theta}{3}+\cos\left(\frac{\theta}{3}+\frac{2\pi}{3}\right)+\cos\left(\frac{\theta}{3}+\frac{4\pi}{3}\right)=0$. By direct computation, we know that \[
    A^3+B^3+C^3-3ABC=(A+B+C)(A^2+B^2+C^2-AB-BC-CA).
\]
This factorization, together with the identity from remark \ref{root},  leads to the following trigonometric identity:
\begin{align*}
	\cos^3\frac{\theta}{3}+\cos^3\left(\frac{\theta}{3}+\frac{2\pi}{3}\right)+\cos^3\left(\frac{\theta}{3}+\frac{4\pi}{3}\right)&=3\cos\frac{\theta}{3}\cos\left(\frac{\theta}{3}+\frac{2\pi}{3}\right)\cos\left(\frac{\theta}{3}+\frac{4\pi}{3}\right)\\
	&=-\frac{3\cos\theta}{4}
\end{align*}
for arbitrary $\theta$.
\end{rem}




\begin{exa}
Find all solutions to $x^3-48x-64\sqrt{2}=0$.
\end{exa}
\begin{sol}
Since $rs=16,~r+s=-4\sqrt{2}$, $r,s$ are the two roots of $t^2+4\sqrt{2}t+16=0$,  we have $r=\frac{-4\sqrt{2}+\sqrt{32-64}}{2}=-2\sqrt{2}+2\sqrt{2}i=4e^{i\frac{3\pi}{4}}$, $s=4e^{-i\frac{3\pi}{4}}$. Hence $\theta=\frac{3\pi}{4}$.  Applying theorem \ref{Case2} we obtain the solutions to the original equations as 
\[
	-2\cdot 4\cos\frac{\pi}{4}=-8\cdot\frac{1}{\sqrt{2}}=-4\sqrt{2},
\]
\[
	-2\cdot 4\cos\left(\frac{\pi}{4}+\frac{2\pi}{3}\right)=8\cos\frac{\pi}{12}=8\cdot\frac{1+\sqrt{3}}{2\sqrt{2}}=2\sqrt{2}+2\sqrt{6},
\]
\[
	-2\cdot 4\cos\left(\frac{\pi}{4}+\frac{4\pi}{3}\right)=-8\cos\frac{19\pi}{12}=-8\cdot \frac{-1+\sqrt{3}}{2\sqrt{2}}=2\sqrt{2}-2\sqrt{6}.
\]
\end{sol}

\begin{exa}~Find all solutions to $x^3-\frac{3}{4}x+\frac{\sqrt{3}}{8}=0$.
\end{exa}
\begin{sol}
Since $rs=\frac{1}{4}$, $r+s=\frac{\sqrt{3}}{2}$, $r,s$ are the two roots of $t^2-\frac{\sqrt{3}}{2}t+\frac{1}{4}=0$, we have  $r=\frac{\sqrt{3}+i}{4}=\frac{1}{2}e^{i\frac{\pi}{6}}$, $s=\frac{\sqrt{3}-i}{4}=\frac{1}{2}e^{-i\frac{\pi}{6}}$.  These are complex conjugates and $\theta=\frac{\pi}{6}$. Applying theorem\ref{Case2}, we find that the three solutions to the original equations are
\[
-2\sqrt{\frac{1}{4}}\cos\left(\frac{1}{3}\cdot\frac{\pi}{6}\right)=-\cos\left(\frac{\pi}{18}\right)=\sin\left(\frac{14\pi}{9}\right)
\]
\[
-\cos\left(\frac{\pi}{18}+\frac{2\pi}{3}\right)=-\cos\left(\frac{13\pi}{18}\right)=\sin\left(\frac{2\pi}{9}\right)
\]
\[
-\cos\left(\frac{\pi}{18}+\frac{4\pi}{3}\right)=-\cos\left(\frac{25\pi}{18}\right)=\sin\left(\frac{8\pi}{9}\right).
\]
\end{sol}

\begin{exa}~Find all solutions to $x^3-\frac{3}{4}x+\frac{1}{8}=0$.
\end{exa}
\begin{sol}
Since $rs=\frac{1}{4}$, $r+s=\frac{1}{2}$, we know that $r,s$ are the roots of $t^2-\frac{1}{2}t+\frac{1}{4}=0$.  Then $r=\frac{1+\sqrt{3}i}{4}=\frac{1}{2}e^{i\frac{\pi}{3}}$, $s=\frac{1-\sqrt{3}i}{4}=\frac{1}{2}e^{-i\frac{\pi}{3}}$ are complex conjugates and $\theta=\frac{\pi}{3}$. Applying theorem\ref{Case2}, we find that the three solutions to the original equations are

\[
-2\sqrt{\frac{1}{4}}\cos\left(\frac{1}{3}\cdot\frac{\pi}{3}\right)=-\cos\left(\frac{\pi}{9}\right)=\cos\left(\frac{8\pi}{9}\right)
\]
\[
-\cos\left(\frac{\pi}{9}+\frac{2\pi}{3}\right)=-\cos\left(\frac{7\pi}{9}\right)=\cos\left(\frac{2\pi}{9}\right)
\]
\[
-\cos\left(\frac{\pi}{9}+\frac{4\pi}{3}\right)=-\cos\left(\frac{13\pi}{9}\right)=\cos\left(\frac{4\pi}{9}\right).
\]
\end{sol}

\medskip

\section{Conclusion}


\medskip

\noindent

Considering the results we have so far, one notices that although we divide the solutions to the equation $x^3-3rsx+rs(r+s)=0$ into three cases ($r$ and $s$ are equal;  $r$ and $s$ are  distinct real numbers; $r$ and $s$ are a pair of complex conjugates), the solutions we obtain in these three cases can actually be expressed uniformly as (\ref{Roots}).  Therefore we may summarize  our results as follows:\\

\medskip

\begin{theorem}
Given a polynomial equation of degree three $x^3+px+q=0$~($p,q\neq 0$) with real coefficients, we can rewrite it as $x^3-3rsx+rs(r+s)=0$.  Then the three solutions to this equation are 
\[
	-r^{\frac{1}{3}}s^{\frac{1}{3}}(r^{\frac{1}{3}}+s^{\frac{1}{3}}),
	~-r^{\frac{1}{3}}s^{\frac{1}{3}}(\omega r^{\frac{1}{3}}+\omega^2 s^{\frac{1}{3}}),
	~-r^{\frac{1}{3}}s^{\frac{1}{3}}(\omega^2 r^{\frac{1}{3}}+\omega s^{\frac{1}{3}}).
\]
In particular, when $r=s$ the three solutions can be further simplified as
\[
	-2\sqrt{rs},\sqrt{rs},\sqrt{rs}.
\]
When $r\neq s$, and $r$ and $s$ are complex conjugates, we can let $\theta={\rm Arg}(r)$.  Then the three solutions can be further simplified as
\[
	-2\sqrt{rs}\cos\frac{\theta}{3},~-2\sqrt{rs}\cos\left(\frac{\theta}{3}+\frac{2\pi}{3}\right),~-2\sqrt{rs}\cos\left(\frac{\theta}{3}+\frac{4\pi}{3}\right).
\]
\end{theorem}

\section{Acknowledgements}
The authors would like to thank Professor William Y.C. Chen for his valuable comments and encouragement, and for offering the newest preprint of \cite{YCChen}. We would also like to thank anonymous referee for helpful comments. 
\rm
\bigskip
\rm

\end{document}